\documentclass[a4paper,11pt]{amsart}
\usepackage{amsmath, amsthm, amssymb}
\usepackage[percent]{overpic}
\usepackage{tikz}
\usepackage{graphicx}
\usepackage{mathbbol}
\usepackage{txfonts}
\usepackage{srcltx} % JD - to jump between source and preview

\newtheorem{thm}{Theorem}[section]
\newcommand{\bt}{\begin{thm}}
\newcommand{\et}{\end{thm}}

\newtheorem*{theorem-non}{Theorem}

\newtheorem{cor}[thm]{Corollary}   %remember switch all {coro} to {cor}
\newcommand{\bc}{\begin{cor}}
\newcommand{\ec}{\end{cor}}

\newtheorem{lem}[thm]{Lemma}   %remember to switch all {lem} to {lem}
\newcommand{\bl}{\begin{lem}}
\newcommand{\el}{\end{lem}}

\newtheorem{prop}[thm]{Proposition}
\newcommand{\bp}{\begin{prop}}
\newcommand{\ep}{\end{prop}}

\newtheorem{defn}[thm]{Definition}
\newcommand{\bd}{\begin{defn}}    % This produces an error????    
\newcommand{\ed}{\end{defn}}

\newtheorem{rmrk}[thm]{Remark}   %remember to switch all {rmrk} to {rmrk}

\newcommand{\br}{\begin{rmrk}}
\newcommand{\er}{\end{rmrk}}

\newcommand{\be}{\begin{equation}}

\newcommand{\ee}{\end{equation}}

\begin{document}

%\doublespacing

\title{A characterization of Whitney forms}

\author{J\'ozef Dodziuk}
\address{Ph.D.~Program in Mathematics\\ CUNY Graduate Center}
\thanks{J\'ozef Dodziuk gratefully acknowledges the hospitality of the Einstein Institute of Mathematics at the Hebrew University of Jerusalem while this note was written. He is grateful to Marta Lewicka for insisting that the theorem deserves a proof.}
\email{jdodziuk@gmail.com}

%\date


%49Q15 (Geometric measure and integration theory, integral and normal currents)


\begin{abstract}
We give a characterization of Whitney forms on an $n$-simplex $\sigma$ and prove that for
every real valued simplicial $k$-cochain $c$ on $\sigma$, the form $Wc$ is the unique differential $k$-form $\varphi$ on $\sigma$ with affine coefficients that pulls back to a constant form of degree $k$ on every
$k$-face $\tau$ of $\sigma$ and satisfies $\int_{\tau} \varphi = <c,\tau>$.

\end{abstract}

\date{\today}

\maketitle

\section{Introduction}
Whitney forms have been extraordinarily useful in 
several areas of mathematics: algebraic topology
\cite{Sullivan}, \cite{Lueck}; global analysis and spectral geometry \cite{Dodziuk1}, \cite{Dodziuk2}; numerical electromagnetism\cite{Bossavit1}, \cite{Bossavit2}; vibrations of thin plates \cite{Simanca}.
Their definition in Whitney's book {\cite[p.~140]{Whitney}
appears somewhat mysterious. Attempts to gain a better insight into the definition have continued up to now. For example, the recent paper of Lohi and Kettunen
\cite{Lohi-Kettunen} contains \emph{three different equivalent definitions.} In this note we give
a conceptual, easily stated characterization of Whitney forms.

On a triangulated differentiable manifold $M$ of $n$ dimensions with
a triangulation $h: K \longrightarrow M$, cf. \cite[p.~124]{Whitney}, the Whitney form $Wc$ corresponding to the
cochain $C^k(K)$ is a family $\omega_\sigma$ of smooth
$k$-forms, satisfying 
certain compatibility conditions, on each closed $n$-simplex $\sigma$. Namely,
if $\tau$ is a common face of two top dimensional
faces $\sigma_1$ and $\sigma_2$, than the pull-backs
to $\tau$ of $\omega_{\sigma_1}$ and $\omega_{\sigma_2}$ coincide. Thus to describe the
Whitney form $Wc$ it suffices to give a description
of $Wc\, |\,\sigma=\omega_\sigma$ for every simplex $\sigma$ of top dimension. Note that the
homeomorphism $h$ defines an affine structure on $\sigma$ and the induced affine structures on common faces of two $n$-simplexes agree. Thus the concept of an affine function on a simplex is well-defined and so is a notion of a "constant" form of degree
$k$ on a $k$-simplex.

From now on we
work on a fixed $n$-simplex $\sigma$. Our characterization of $Wc$ is stated precisely in the
Theorem below. It asserts that $Wc$ restricted to $\sigma$ is the unique $k$-form on $\sigma$ with affine coefficients and constant pull-backs to $k$-faces whose integrals over $k$-faces $\tau$ are prescribed by the values $\langle c,\tau \rangle$ of $c$ on $\tau$.

\section{Proof of the Theorem}
A simplex $\tau=[p_0,p_1, \ldots ,p_k]$ of $k$ dimensions is a convex hull of $k+1$ points in general position in $\mathbb{R}^n$. In particular, every simplex is closed. We will consider a fixed $n$-simplex $\sigma$ together with all its $k$-faces $\tau$ with $0\leq k \leq n$. Thus a point $q\in \sigma$ is a
convex linear combination
\begin{gather*}
q = m_0p_0 + m_1p_1 +\ldots + m_n p_n\\
m_i  \geq 0\quad \mbox{for}\quad i=0,1\ldots ,n\\
 m_0+ m_1 \ldots  +m_n =1
\end{gather*}
and the barycentric coordinate functions $\nu_i(q)$ are defined by
$$
\nu_i(q) = m_i .
$$
We observe that, if $q=(x^1,x^2, \ldots , x^n)$ the barycentric coordinates are affine functions of $x^1, x^2, \ldots , x^n$ i.e.\ are of the form 
$a_1 x^1 +a_2 x^2 + \ldots + a_n  x^n + b$.
We regard all simplices as oriented with the orientation determined by the order of
vertices with the usual convention that $-\tau$ is $\tau$ with the opposite orientation and
that under a permutation of vertices the orientation changes by the sign of the permutation.
A cochain $c$ of degree $k$ is then defined as a formal linear combination with real
coefficients of duals the $\tau^*$ of $k$-faces $\tau$ of $\sigma$ and we denote by $C^k(\sigma) = C^k$ the space of all such cochains. If $c=\sum_{\tau} a_\tau \tau^*$ we will write $a_\tau =\langle c,\tau\rangle$. Finally, we will denote by
$\Lambda^k(\sigma)=\Lambda^k$ the space of all
smooth exterior differential forms of degree $k$ on the simplex $\sigma$.
With this notation, one defines the Whitney mapping
$$W: C^k \longrightarrow \Lambda^k$$
for all $k=0,1,\ldots n$, cf.\ \cite{Whitney} or \cite{Dodziuk2} for a detailed discussion. We will call forms in the image of $W$ the Whitney forms. It follows immediately from the 
definition that the Whitney forms when expressed in terms of the coordinates of $\mathbb{R}^n$ have affine coefficients. We abuse the language and say that a form $\eta\in \Lambda^k(\tau )$  is constant if it is a constant multiple of the
Euclidean volume element on $\tau$. After these preliminaries we state our theorem.

\begin{theorem-non}
Let $\sigma$ be a simplex of $n$ dimensions and $c$ a cochain of degree $k$ on $\sigma$. $Wc$ is the unique $k$-form $\omega$ on $\sigma$
satisfying the following conditions.
\begin{enumerate}
\item
$\omega$ has affine coefficients.
\item The pull-back $\iota^*_\tau\omega$ is constant for every $k$-dimensional face $\tau$ of 
$\sigma$, where $\iota_\tau : \tau \hookrightarrow \sigma$ denotes the inclusion map.
\item $\int_\tau \omega = \langle c, \tau\rangle$ for every $k$-face $\tau$ of $\sigma$. 
\end{enumerate}
\end{theorem-non}
\begin{proof}
We first observe that without any loss of generality we can assume that $\sigma$ is the
standard simplex in $\mathbb{R}^n$ i.e. is given by
$$
\sigma = \left\lbrace \, (x^1,x^2, \ldots , x^n ) \in \mathbb{R}^n \mid x^i \geq 0\quad \mbox{for}\quad i=1,2, \ldots n\,;\quad \sum^n_{i=0} x^i \leq 1\, \right\rbrace.
$$
Thus $\sigma=[0,e_1,e_2\ldots , e_n]$ where $e_i$ is the point on the $i$-th coordinate axis
with $x^i=1$.
The barycentric coordinate functions restricted to $\sigma$ are then given by 
\begin{equation}\label{bar-coord}
\nu_0  = 1-(x^1 + x^2 + \ldots +x^n)\quad\mbox{and}\quad
\nu_i  = x^i \quad \mbox{for} \quad  i=1,2, \ldots n.
\end{equation}

We first do a quick dimension count that makes the theorem plausible.
The dimension of the space of $k$-forms with affine coefficients on $\sigma$ is $
\binom{n}{k}(n+1)$. Requiring that $\iota_\tau^*\omega$ is \emph{constant} on a $k$-simplex $\tau$ imposes $k$ conditions and
the number of $k$-faces of an $n$-simplex is $\binom{n+1}{k+1}$. Thus, the dimension of the space of $k$-forms
satisfying (1) and (2) above ought to be
$$
\binom{n}{k}(n+1) - \binom{n+1}{k+1} k= \binom{n+1}{k+1}.
$$
This last integer is the number of $k$-faces of $\sigma$, i.e.\  the dimension of the space  $C^k(\sigma)$ of $k$-cochains.

It is  instructive to consider the simplest cases $k=0$ and $k=n$ of the theorem.
A $0$-cochain is a sum $c=\sum a_i p_i^*$ and 
\begin{align*}
Wc &=a_0\nu_0 + a_1 \nu_1 +\ldots a_n \nu_n\\
   &= a_0 \left (1-\sum_{i=1}^n x^i \right) + \sum_{i=1}^n a_i x^i\\
   &= a_0 + \sum_{i=1}^n (a_i -a_0) x^i
\end{align*} 
is the unique affine function $f$ taking prescribed values
$f(p_i)=\int_{p_i} f = \langle  c,p_i \rangle$, where the integration of
a form of degree 0 over a vertex is just the evaluation.

If $k=n$, $\sigma$ is the only face of dimension $n$ so every cochain is a multiple 
of $\sigma^*$. For $c=\sigma^*$, we have
\begin{align*}
Wc &= W\sigma^*\\
  &=\left (
  n! \sum_{j=0}^n (-1)^j \nu_j d\nu_0\wedge \ldots \wedge \widehat{d\nu_j} \wedge \ldots         
        \wedge d\nu_n 
        \right )\\
  &=n!  dx^1\wedge \ldots \wedge dx^n
\end{align*}
where we used the explicit expressions of the barycentric coordinates (\ref{bar-coord}) in terms of the coordinates $x^1, \ldots , x^n$ and the hat over a factor means that the factor is omitted. Since the volume of the standard $n$-simplex 
in $\mathbb{R}^n$ is equal to $1/n!$, $\int_\sigma W(\sigma^*) = \langle \sigma^*,\sigma \rangle =1$, $W\sigma^*$ is the unique constant form with prescribed integral equal to one.

We now consider the case when $1 \leq k \leq n-1$. We will write $\Lambda_e^k$ for the space of $k$-forms on $\sigma$ with affine coefficients and with constant pull-backs to $k$-faces of $\sigma$. It is obvious from the definition of $Wc$ and from (\ref{bar-coord})
that $Wc$ has affine coefficients on $\sigma$ for every $c\in C^k(\sigma)$. Similarly,
since $\iota_\tau^* W(c)$ is a form of maximal degree on $\tau$, the calculation above,
with $k$ replacing $n$, shows that $\iota_\tau^* W(c)$ is constant on $\tau$ for every 
$k$-face $\tau$ of $\sigma$. It follows that $WC^k \subset \Lambda^k_e$. Now let $\varphi\in \Lambda^k_e $. We use the restriction of the de Rham map 
$R:\Lambda^k(\sigma) \longrightarrow C^k(\sigma)$,
$$
\langle R\omega , \tau \rangle = \int_\tau \omega,
$$
to $\Lambda^k_e$ and consider the difference $\eta = \varphi - WR \varphi$. Clearly,
$\eta \in \Lambda^k_e$.  Moreover basic properties of the Whitney mapping
(cf. \cite{Whitney, Dodziuk2}) imply that $R\eta = R\varphi - RWR\varphi=R\varphi-R\varphi=0$, i.e. $\eta$ integrates to zero on every $k$-face
of $\sigma$. Since the pull-back $\iota^*\eta$ is constant on every such face $\tau$, $\iota^*_\tau \eta$
vanishes identically on every $k$-face $\tau$. Thus to show that $\varphi = WR \varphi$ (which
would prove our theorem) it suffices to show that every form $\eta\in \Lambda^k_e$, whose
pull-backs to all $k$-faces vanish, is itself identically zero on $\sigma$. Let $\eta$ be
such a form. We express it in the standard coordinates of $\mathbb{R}^n$ as follows.
\begin{equation}
\label{eta:coord}
\eta =  \sum_I( b_I +a_{I,1}x^1 + \ldots +a_{I,n}x^n )dx^I
\end{equation}
Here $I$ is a multi-index $I=(i_1<i_2 < \ldots <i_k)$,
 $1\leq i_j \leq n$ for every $j$ and
$dx^I = dx^{i_1}\wedge dx^{i_2} \wedge \ldots\wedge dx^{i_k} $. We will abuse the notation at times and think of $I$ as a set. Fix a multi-index $J$ and consider the coordinate plane of the variables $x^{j_1},x^{j_2},\ldots , x^{j_k}$. 

Let $\tau_J$ denote the $k$-face of $\sigma$ contained in that plane. By assumption $\iota^*_{\tau_J}\eta$ is identically zero. The variables $x_t$ for $t\not\in J$ vanish in
this plane so that
\begin{equation}
\label{eta:coord:plane}
\iota^*_{\tau_J} \eta\, = \sum_{t\in J} (a_{J,t}x^t + b_J) dx^J \equiv 0.
\end{equation}
Since $J$ was arbitrary, $b_J=0$ and $a_{J,t}=0$ for all $J$ and all $t\in J$. It follows that we can
rewrite (\ref{eta:coord}) on $\sigma$ as follows.
\begin{equation}
\label{eta:inclined}
\eta =  \sum_I\sum_{j\not\in I}a_{I,j}x^j  dx^I
\end{equation}
Again, fix the multi-index $L$, an integer $m\not\in L$, $1\leq m \leq n-1$, and the simplex $\tau=[e_m,e_{l_1},\ldots ,e_{l_k}]$. $\tau$ is a $k$-simplex in the $(k+1)$-plane $P$ with coordinates $x^m,x^{l_1}, \ldots ,x^{l_k}$ as in the figure below.
Recall that on $\tau$, $x^{l_1}, \ldots ,x^{l_k}$ can be taken as local
coordinates since 
\begin{equation}\label{xm}
x^m = 1 - (x^{l_1} + \ldots + x^{l_k})
\end{equation}
Moreover
\begin{equation}\label{dxm}
dx^m= -(dx^{l_1} + \ldots + dx^{l_k})
\end{equation}
\begin{figure}[h!]
\vspace*{1ex}
\usetikzlibrary {arrows.meta}
\tikzset{help lines/.style=thick}
\begin{tikzpicture}[>=Stealth]
  \draw[->] (0,0) -- (5,0);
  \draw[->] (0,0) -- (0,3.5);
  \draw[->]  (0,0) -- (-1,-1);
  \draw (0,3) -- (3,0);
  \draw (0.5,3.2) node {$x_m$};
  \draw (1.7,1.7) node {$\tau$};
  \draw (4.5,0.3) node {$x_{l_1},x_{l_2},\ldots, x_{l_k}$};
  \draw (0.69,-0.8) node {$x_j,j\not\in L\cup \lbrace m \rbrace$};
  \draw (2.5,2.5) node {$P$};
\end{tikzpicture}
\end{figure}
We express the pull-back $\iota_\tau^* \eta$ in terms these coordinates using (\ref{xm}) and
(\ref{dxm}). Observe that if
$I\cup \lbrace j \rbrace \neq L \cup \lbrace m \rbrace$ 
one of the indices in $I\cup \lbrace j \rbrace$ is not in $L \cup \lbrace m \rbrace$. 
The corresponding variable is identically zero on the plane $P$ so that the
summand $a_{I,j}x^jdx^I$ vanishes on $P$ and is therefore equal to zero when pulled back to
$\tau$. 
Therefore 
\begin{equation}\label{pullback:tau}
\iota^*_\tau \eta \,\, = \sum_{I\cup \lbrace j \rbrace = L \cup \lbrace m \rbrace}
a_{I,j}x^jdx^I.
\end{equation}
Now consider the summand with $I=L$ and $j=m$. The coefficient of $dx^L$ in this term
is
$$
a_{L,m}x^m+a_{L,l_1}x^{l_1}+\ldots +a_{L,l_k}x^{l_k}
$$
and we use (\ref{xm}) to eliminate $x^m$.

Thus, on $\tau$, the coefficient in question can be written as 
$$
a_{L,m} - a_{L,m}\sum_{s=1}^k x^{l_s} + a_{L,l_1}x^{l_1}+\ldots +a_{L,l_k}x^{l_k}.
$$
Remaining terms in the sum (\ref{pullback:tau}) have $j\neq m$. It follows that, for
those terms, $x^j$ is one of $x^{l_1}, \ldots ,x^{l_k}$ and $x^m$ enters only into
the differential monomial $dx^I$ from which it can be eliminated using (\ref{dxm}).
It follows that
$$
\iota^*_\tau \eta = (\,a_{L,m} + \mbox{\it linear terms}\, )\, dx^L.
$$
Since $\iota^*_\tau \eta$ is assumed to be identically zero, $a_{L,m}=0$. $L$ was fixed but arbitrary so that
$\eta \equiv 0$.

\end{proof}

\bibliographystyle{plain}

\bibliography{whitney-characterization}

\begin{thebibliography}{1}

\bibitem{Bossavit1}
A.~Bossavit.
\newblock A uniform rationale for {W}hitney forms on various supporting shapes.
\newblock {\em Math. Comput. Simulation}, 80(8):1567--1577, 2010.

\bibitem{Bossavit2}
Alain Bossavit.
\newblock {\em Computational electromagnetism}.
\newblock Electromagnetism. Academic Press, Inc., San Diego, CA, 1998.
\newblock Variational formulations, complementarity, edge elements.

\bibitem{Dodziuk2}
Jozef Dodziuk.
\newblock Finite-difference approach to the {H}odge theory of harmonic forms.
\newblock {\em Amer. J. Math.}, 98(1):79--104, 1976.

\bibitem{Dodziuk1}
Jozef Dodziuk.
\newblock de {R}ham-{H}odge theory for {$L^{2}$}-cohomology of infinite
  coverings.
\newblock {\em Topology}, 16(2):157--165, 1977.

\bibitem{Lohi-Kettunen}
Jonni Lohi and Lauri Kettunen.
\newblock Whitney forms and their extensions.
\newblock {\em J. Comput. Appl. Math.}, 393:Paper No. 113520, 19, 2021.

\bibitem{Lueck}
Wolfgang L\"{u}ck.
\newblock {\em {$L^2$}-invariants: theory and applications to geometry and
  {$K$}-theory}, volume~44 of {\em Ergebnisse der Mathematik und ihrer
  Grenzgebiete. 3. Folge. A Series of Modern Surveys in Mathematics [Results in
  Mathematics and Related Areas. 3rd Series. A Series of Modern Surveys in
  Mathematics]}.
\newblock Springer-Verlag, Berlin, 2002.

\bibitem{Simanca}
Santiago~R. Simanca.
\newblock The (small) vibrations of thin plates.
\newblock {\em Nonlinearity}, 32(4):1175--1205, 2019.

\bibitem{Sullivan}
Dennis Sullivan.
\newblock Cartan-de {R}ham homotopy theory.
\newblock In {\em Colloque ``{A}nalyse et {T}opologie'' en l'{H}onneur de
  {H}enri {C}artan ({O}rsay, 1974)}, pages 227--254. Ast\'{e}risque, No.
  32--33. 1976.

\bibitem{Whitney}
Hassler Whitney.
\newblock {\em Geometric integration theory}.
\newblock Princeton University Press, Princeton, N. J., 1957.

\end{thebibliography}

\end{document}